\documentclass[11pt,a4paper]{article}

\usepackage{latexsym}
\usepackage{mathrsfs}
\usepackage{extarrows}
\usepackage{graphicx}
\usepackage{supertabular}
\usepackage{subfigure}
\usepackage{setspace}
\usepackage{epstopdf}
\usepackage{amssymb}
\usepackage{mathrsfs}
\usepackage{amsmath}
\usepackage{amsfonts, amsthm, amssymb}
\usepackage{amsfonts}
\usepackage{color}
\usepackage{float}
\usepackage{graphicx}
\usepackage{subfigure}
\usepackage{caption}
\usepackage{amsmath}
\numberwithin{equation}{section}
\newtheorem{lem}{Lemma}[section]
\newtheorem{thm}[lem]{Theorem}
\newtheorem{cor}[lem]{Corollary}
\newtheorem{ques}[lem]{Question}
\newtheorem{claim}{Claim}[section]

\newtheorem*{rem}{Remark}
\newtheorem{fact}[lem]{Fact}
\newtheorem{defi}[lem]{Definition}

\begin{document}

\title{Chromatic profiles of odd cycles}
\author{ Zilong Yan \footnote{School of Mathematics, Hunan University,
 Changsha 410082, P. R. China. E-mail: zilongyan@hnu.edu.cn. Supported in part by Postdoctoral fellow fund in China (No. 2023M741131)  and  Postdoctoral Fellowship Program of
CPSF under Grant (No. GZC20240455).} ~~~~ Yuejian Peng \footnote{Corresponding author. School of Mathematics, Hunan University, Changsha 410082, P. R. China. E-mail: ypeng1@hnu.edu.cn. Supported in part by National Natural Science Foundation of China (No. 11931002 and 12371327)}  ~~~~  Xiaoli Yuan \footnote{School of Mathematics, Hunan University,
 Changsha 410082, P. R. China. E-mail: xiaoliyuan@hnu.edu.cn.} \\
}
\date{}
\maketitle
\begin{abstract}
Erd\H{o}s and Simonovits  asked the following question: For an integer $c\geq 2$ and a family of non-bipartite graphs $\mathcal{F}$, what is the infimum of $\alpha$ such that any $\mathcal{F}$-free $n$-vertex graph with $n$ large enough and  minimum degree at least $\alpha n$ has chromatic number at most $c$? Denote the infimum as $\delta_{\chi}(\mathcal{F}, c)$.  A fundamental  result of
 Erd\H{o}s, Stone and Simonovits  implies that if $3\le r+1=\chi(\mathcal{F})=\min\{\chi (F): F\in \mathcal{F}\}$, then
  for any $c\le r-1$,
$\delta_{\chi}(\mathcal{F}, c)=1-{1 \over r}$.  So the remaining challenge is to determine $\delta_{\chi}(\mathcal{F}, c)$ for $c\ge \chi (\mathcal{F})-1$. Most previous known results are under the condition that $c= \chi (\mathcal{F})-1$.  When $c\ge \chi (\mathcal{F})$, the only known exact results are  $\delta_{\chi}(K_3, 3)$ by H\"aggkvist and Jin, and $\delta_{\chi}(K_3, c)$ for every $c\ge4$ by Brandt and Thomass\'{e}, $\delta_{\chi}(K_r, r)$ and $\delta_{\chi}(K_r, r+1)$ by Goddard and Lyle,  and Nikiforov.  Combining some results of Thomassen and Ma,  $\Omega\bigg((c+1)^{-8(k+1)}\bigg)=\delta_{\chi}(C_{2k+1}, c)=O(\frac{k}{c})$ for $c\ge 3$.  In this paper, we determine $\delta_{\chi}(C_{2k+1}, c)$ for all $c\ge 2$ and $k\ge 3c+4$ ($k\ge 5$ if $c=2$). We also obtain the following corollary. If $G$ is a  graph on $n$ vertices with $c\ge 3$, $\chi(G)>c$ and $\delta(G)> {n \over 2c+2}$, then $C_{2k+1} \subset G$ for all $k\in [3c+4, {n \over 108(c+1)^c}]$. Methods to obtain all previous known results related to odd cycles cannot be applied to solve for $\delta_{\chi}(C_{2k+1}, c)$ for $c\ge 3$.
The innovation  of our proof is to give the concept of a `strong $2k$-core'.  We think that this concept grasps the essence of the problem and it makes  our proof  concise and elementary (we do not need to borrow any other tools). How to define a proper `core' might be a key to this type of questions.

\end{abstract}
\section{Introduction}
For  a graph $G$, let $V(G)$ and $E(G)$ denote the vertex set and the edge set of $G$ respectively.
Denote the minimum degree,  the chromatic number of a graph $G$  and the clique number of $G$ as $\delta(G)$,   $\chi(G)$ and $\omega(G)$ respectively. Let $\mathcal{F}$ be a family of  graphs, a  graph $G$ is {\em $\mathcal{F}$-free} if $G$ does not contain any member of $\mathcal{F}$ as a subgraph. If $\mathcal F$ consists of a single graph $F$, we simply write  $\{F\}$-free as $F$-free. Let $K_{r+1}$ denote a complete graph on $r+1$ vertices, and let $C_{k}$ denote a cycle with $k$ vertices.

Andr\'{a}sfai, Erd\H{o}s and S\'{o}s \cite{Andrsfi}  showed that if $G$ is a $K_{r+1}$-free graph on $n$ vertices with  $\delta(G)>\frac{3r-4}{3r-1}n$, then  $\chi(G)\le r$. This result is very interesting  since the difference between the chromatic number of a graph and its clique number could be very large, and this result tells us that $\omega(G)\le r$ implying that
$\chi(G)\le r$  if $G$ has large enough minimum degree.
Inspired by the work of  Andr\'{a}sfai, Erd\H{o}s and S\'{o}s \cite{Andrsfi}, Erd\H os and Simonovits \cite{Erdos}  asked the following  general question.

\begin{ques}[Erd\H os and Simonovits \cite{Erdos}]
For a family $\mathcal{F}$ of  graphs and a positive integer $c$, what is the minimum  $f$ such that an $\mathcal{F}$-free graph $G$ on $n$ vertices with $\delta(G)\ge f$ must satisfy $\chi(G)\le c$?
\end{ques}

This question is interesting in several perspectives. If we take  $\mathcal{F}=\{K_{r+1}\}$ as in  \cite{Andrsfi}, then the question asks to determine a tight minimum degree condition on a graph $G$ such that  $\chi(G)\le c$ if $\omega(G)\le r$.
Erd\H os showed that for any $k$ and $c$,  there is a graph $G$ such that the length of a shortest cycle in $G$ is at least $k+1$ and $\chi(G)>c$. If we take $\mathcal{F}=\{C_3, C_4, \cdots, C_k\}$, then the question of Erd\H os and Simonovits \cite{Erdos} asks to determine  a tight minimum degree condition on a graph $G$ such that the length of a shortest cycle in $G$ is at least $k+1$ but $\chi(G)\le c$.
  An equivalent form of the question  of  Erd\H os and Simonovits \cite{Erdos} is:
 for a family of  graphs $\mathcal{F}$ and a positive integer $c$, what is the minimum  $f$ such that a graph $G$ on $n$ vertices with  $\delta(G)\ge f$ and $\chi(G)>c$ must contain a copy of some graph in  $\mathcal{F}$?
A classical question in graph theory is to determine a tight minimum degree condition to guarantee the existence of certain subgraphs.
For example, the existence of cycles in graphs with sufficiently large minimum degree has been intensively studied (see \cite{alon,Bal,Brandt,Er,Hggkvist,Illingworth,KSV,Let,Nikiforov,Sankar,Ver}). The question of Erd\H os and Simonovits \cite{Erdos}  asks what happens if we know some information on the chromatic number $\chi(G)$.

The question of  Erd\H os and Simonovits \cite{Erdos} has been stimulating the general study of the so-called chromatic profile, and it also has connections with the stability of Tur\'{a}n problems. Let us give details below.
Denote $$\mathcal G(n, \mathcal F, \alpha)=\{G: G \ \text{is} \ \mathcal F\text{-free} \ \text{with} \ n \ \text{vertices  and } \  \delta(G)\ge \alpha n. \}.$$
Call $\mathcal{F}$ {\em non-degenerate} if $\chi(F)>2$ for each $F\in \mathcal{F}$.
\begin{defi}
For  a positive integer $c\ge 2$, the  chromatic profile of a non-degenerate family $\mathcal F$ as a function in $c$ is defined to be
$$\delta_{\chi}(\mathcal{F}, c)=\inf\{\alpha| \ \text{Any} \ G\in \mathcal G(n, \mathcal F, \alpha)  \ \text{with} \ n \   \text{ large enough must satisfy} \ \chi(G)\le c.\}.$$
If $\mathcal F$ consists of a single graph $F$, we simply write  $\delta_{\chi}(\{F\}, c)$ as $\delta_{\chi}(F, c)$.
\end{defi}

 The study of chromatic profile is also related to the stability of Tur\'{a}n problems. For  a family $\mathcal{F}$ of $k$-uniform graphs and a positive integer $n$, the {\em Tur\'{a}n number} $ex(n, \mathcal{F})$ is the maximum number of edges  an $\mathcal{F}$-free $k$-uniform graph on $n$ vertices can have.    An averaging argument of Katona, Nemetz and Simonovits \cite{KNS} shows that the sequence ${ex(n,\mathcal F) \over {n \choose k }}$ is  non-increasing. Hence $\underset{n\rightarrow\infty}{\lim}{ex(n, \mathcal F) \over {n \choose k } }$ exists. The {\em Tur\'{a}n density} of $\mathcal F$ is defined as $\pi(\mathcal F)=\underset{n\rightarrow\infty}{\lim}{ex(n, \mathcal F) \over {n \choose k } }.$
If $\mathcal F$ consists of a single $k$-uniform graph $F$, we simply write $ex(n, \{F\})$ and $\pi(\{F\})$ as $ex(n, F)$ and $\pi(F)$.
   Let $T_r(n)$ denote the complete $r$-partite graph on $n$ vertices where
 its part sizes are as equal as possible.
Tur\'{a}n \cite{Tur}   obtained that if $G$ is a $K_{r+1}$-free graph on $n$ vertices,
then $e(G)\le e(T_r(n))$, equality holds if and only if $G=T_r(n)$.
  A fundamental  result of
 Erd\H{o}s, Stone and Simonovits \cite{ES46, ES66} gives the asymptotical value of $ex(n, \mathcal{F})$ for all non-degenerate families of graphs. For a family $\mathcal{F}$ of graphs, let $\chi(\mathcal{F})=\min\{\chi (F): F\in \mathcal{F}\}$.

\begin{thm}[Erd\H{o}s--Stone--Simonovits \cite{ES46, ES66}]  \label{thm11}
If $\mathcal{F}$ is a family of graphs with
$\chi (\mathcal{F})=r+1$, then
\[  \mathrm{ex}(n,\mathcal{F}) =e(T_r(n)) + o(n^2)= \left(  1-\frac{1}{r} + o(1)\right) \frac{n^2}{2}. \]
\end{thm}
Erd\H{o}s \cite{Erd1966Sta1,Erds} and Simonovits \cite{Sim1966} also obtained
a stronger structural  theorem of Theorem \ref{thm11} and
discovered  a certain stability phenomenon.

\begin{thm}[Erd\H{o}s--Simonovits \cite{Erd1966Sta1,Erds, Sim1966}] \label{thm12}
Let $\mathcal{F}$ be a family of graphs with
$\chi (\mathcal{F})=r+1\ge 3$.
For every $\varepsilon >0$,
there exist $\delta >0$ and $n_0$ such that
if $G$ is a graph on $n\ge n_0$ vertices,
 and $G$ is $\mathcal{F}$-free such that
$e(G)\ge (1- \frac{1}{r} - \delta) \frac{n^2}{2}$,
 then  $G$ can be obtained from $T_r(n)$
by adding or deleting a total number of at most $\varepsilon n^2$ edges.
\end{thm}

Recently, F\"uredi\cite{Furedi},  Roberts and Scott \cite{RS21}, Balogh, Clemen, Lavrov, Lidicky and
Pfender \cite{ Balogh},  Kor\'{a}ndi, Roberts and Scott \cite{KRS21} give further developments on the relationship between $\varepsilon$ and $\delta$ in Theorem \ref{thm12}. We do not give statements of their results here, and  refer  readers to their papers.

 Liu,  Mubayi and Reiher \cite{LMR} also defined the degree-stability of Tur\'{a}n problems for $k$-uniform graphs which implies the edge-stability as in Theorem \ref{thm12}.
Let $\mathcal F$ be a non-degenerate family of $k$-uniform graphs (i.e. $\pi(\mathcal F)>0$). Let $\mathscr K$ be a class of $\mathcal F$-free $k$-uniform graphs. If there exists $\epsilon >0$ and $n_0$ such that every $\mathcal F$-free $k$-uniform graph $G$ on $n\ge n_0$ vertices with $\delta(G)\ge (\pi(\mathcal F)-\epsilon){n \choose k-1}$ is a subgraph of some member of $\mathscr K$, then we say that $\mathcal F$ is {\em degree-stable} with respect to $\mathscr K$. It is not hard to see that the degree-stability implies the edge-stability (as in Theorem \ref{thm12}) since we can delete vertices with `small' degrees by losing $o(n^{k-1})$ edges. We can define the chromatic profile for a family of $k$-uniform graphs similarly. Following from the definitions, we can see that for a non-degenerate family $\mathcal{F}$  of $k$-uniform graphs ($\pi(\mathcal F)>0$), $\mathcal F$ is degree-stable with respect to the family of all $r$-colorable $k$-uniform graphs if and only if $\delta_{\chi}(\mathcal{F}, r)\le  \pi(\mathcal F)$.

\begin{rem}\label{csmall}
If $\mathcal{F}$ is a family of graphs with
$\chi (\mathcal{F})=r+1\ge 3$,  then for any $c\le r-1$,
$\delta_{\chi}(\mathcal{F}, c)=1-{1 \over r}$.
\end{rem}
\noindent\emph{\textbf{Proof}.}
Let $\mathcal{F}$ is a family of graphs with
$\chi (\mathcal{F})=r+1\ge 3$. By Theorem \ref{thm11}, any $\mathcal{F}$-free graph with $n$ vertices has minimum degree at most $(1- \frac{1}{r} +o(1))n$. Thus, $\delta_{\chi}(\mathcal{F}, c)\le 1-{1 \over r}$.
On the other hand, the Tur\'{a}n graph $T_r(n)$ is $\mathcal{F}$-free with chromatic number $r$ and $\delta(T_r(n))=(1- \frac{1}{r} +o(1))n$. So $\delta_{\chi}(\mathcal{F}, c)\le 1-{1 \over r}$. This completes the proof. \qed

By Remark \ref{csmall}, the remaining challenge is to determine $\delta_{\chi}(\mathcal{F}, c)$ for $c\ge \chi (\mathcal{F})-1$.
The first study of chromatic profile is due to  Andr\'{a}sfai, Erd\H{o}s and S\'{o}s \cite{Andrsfi}, 
they showed that $\delta_{\chi}(K_{r+1}, r)=\frac{3r-4}{3r-1}$.
In the same paper, they also showed that $\delta_{\chi}(\{C_3,C_5,\cdots, C_{2k+1}\}, 2)=\frac{2}{2k+3}$. H\"{a}ggkvist \cite{Hggkvist} showed that $\delta_{\chi}(C_{2k+1}, 2)=\frac{2}{2k+3}$ for $k=\{1,2,3,4\}$. Recently, Yuan-Peng \cite{YuanPeng, YuanPeng2} unified the results of Andr\'{a}sfai, Erd\H{o}s and S\'{o}s \cite{Andrsfi} and H\"{a}ggkvist \cite{Hggkvist}, we obtained the chromatic profile for any family consisting of some odd cycles, we proved that for a family ${\mathcal C}$ of odd cycles in which $C_{2p+1}$ is the shortest odd cycle not in ${\mathcal C}$, and $C_{2k+1}$ is the longest odd cycle in ${\mathcal C}$,   $\delta_{\chi}({\mathcal C}, 2)=\max\{ \frac{1}{2(2p+1)}, \frac{2}{2k+3}\}$. The lower bound is evidenced by the following  two graphs. Let $BC_{2p+1}(n)$ denote the graph obtained by taking $2p+1$ vertex-disjoint copies of $K_{\frac{n}{2(2p+1)},\frac{n}{2(2p+1)}}$ and selecting a vertex in each of them such that these vertices form a cycle of length $2p+1$. Let $C_{2k+3}({n \over 2k+3})$ denote the balanced blow up of $C_{2k+3}$ with $n$ vertices. See Fig.1 and Fig.2. Note that both $BC_{2p+1}(n)$ and $C_{2k+3}({n \over 2k+3})$ are ${\mathcal C}$-free and non-bipartite, thus, $\delta_{\chi}({\mathcal C}, 2)\ge \max\{ \frac{1}{2(2p+1)}, \frac{2}{2k+3}\}$. In \cite{YuanPeng, YuanPeng2}, we  proved that $\delta_{\chi}({\mathcal C}, 2)\ge \max\{ \frac{1}{2(2p+1)}, \frac{2}{2k+3}\}$.
\begin{figure}[ht]
\centering
\begin{minipage}[b]{0.44\linewidth}
\includegraphics[height=4.5cm]{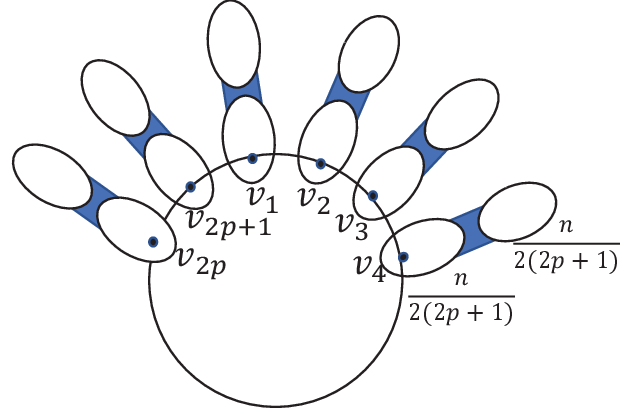}
\caption{$BC_{2p+1}(n)$}
\label{fig1}
\end{minipage}
\quad
\quad
\quad
\begin{minipage}[b]{0.42\linewidth}
\includegraphics[height=4.8cm]{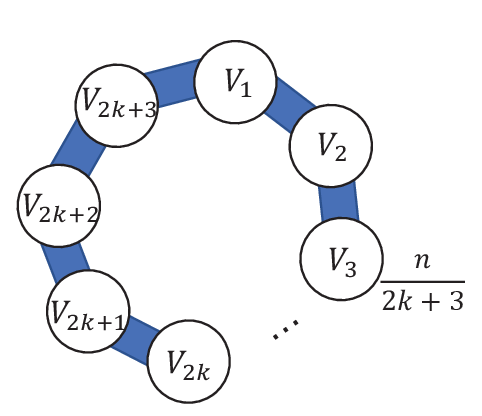}
\caption{$C_{2k+3}({n \over 2k+3})$}
\label{fig2}
\end{minipage}
\end{figure}

Precisely, we \cite{YuanPeng, YuanPeng2}  showed the following result.
\begin{thm}[Yuan-Peng \cite{YuanPeng, YuanPeng2}] \label{YuanPeng}
Let  $n\ge 1000k^{8}$ be positive integers.  Let ${\mathcal C}$ be a family of some odd cycles in which $C_{2p+1}$ is the shortest odd cycle not in ${\mathcal C}$ and $C_{2k+1}$ is the longest odd cycle in ${\mathcal C}$.   If $G$ is an $n$-vertex ${\mathcal C}$-free graph with $\delta(G)>\max\{ \frac{n}{2(2p+1)}, \frac{2}{2k+3}n\}$, then $G$ is bipartite. Furthermore, if $k\ge 4p+1$, then the only $n$-vertex ${\mathcal C}$-free  non-bipartite graph with minimum degree $\max\{ \frac{n}{2(2p+1)}, \frac{2}{2k+3}n\}=\frac{n}{2(2p+1)}$ is $BC_{2p+1}(n)$, and if $k\le 4p$, then the  only $n$-vertex ${\mathcal C}$-free  non-bipartite graph with minimum degree $\max\{ \frac{n}{2(2p+1)}, \frac{2}{2k+3}n\}=\frac{2}{2k+3}n$ is $C_{2k+3}({n \over 2k+3})$.
\end{thm}

When $c\ge \chi (\mathcal{F})$, we know very few results on $\delta_{\chi}(\mathcal{F}, c)$. The only known exact results are  $\delta_{\chi}(K_3, 3)=\frac{10}{29}$ by H\"aggkvist \cite{Hggkvist} and Jin \cite{Jin}, and $\delta_{\chi}(K_3, c)=\frac{1}{3}$ for every $c\ge4$ by Brandt and Thomass\'{e} \cite{BrTh},  $\delta_{\chi}(K_{r+1}, r+1)=1-\frac{19}{19r-9}$ and $\delta_{\chi}(K_{r+1}, r+2)=1-\frac{2}{2r-1}$ by Goddard and Lyle \cite{GoLy}, and Nikiforov \cite{Nikiforov1}.

For odd cycles, we do not know any exact results when $c\ge 3$ except $\delta_{\chi}(K_3, c)$.
Van Ngoc and Tuza \cite{Van} showed that $\delta_{\chi}(\{C_3,C_5,\cdots, C_{2k-1}\},3)\ge \frac{3}{2k^2+k+1}$.
Recently, B\"{o}ttcher, Frankl, Cecchelli, Parzcyk and Skokan \cite{Bottcher} showed the following upper bound.
\begin{thm}[B\"{o}ttcher, Frankl, Cecchelli, Parzcyk and Skokan \cite{Bottcher}]
Let $t$ be a positive integer and $k\ge 45t+5490$. Then $\delta_{\chi}(\{C_3,C_5,\cdots, C_{2k-1}\},
3)\le \frac{1}{2k+t}$.
\end{thm}
Thomassen \cite{Thomassen} showed that $\delta_{\chi}(C_5, c)\le\frac{6}{c}$ and gave an upper bound for $\delta_{\chi}(C_{2k+1}, c)$, together with a result of Ma \cite{Ma} (it is not the main focus in \cite{Ma}), we have
\begin{rem}(Thomassen \cite{Thomassen}, Ma \cite{Ma})
$$\Omega\bigg((c+1)^{-8(k+1)}\bigg)=\delta_{\chi}(C_{2k+1}, c)=O(\frac{k}{c}).$$
\end{rem}

\section{Main results and essence of the proof}

In this paper we show that 
$\delta_{\chi}(C_{2k+1}, r)=\frac{1}{2r+2}$ for all $r\ge 3$ and $k\ge 3r+4$.
 We also improve the requirement on $n$ in Theorem \ref{YuanPeng} when $k\ge 4p+1$.
Let us state our main results precisely.

To determine $\delta_{\chi}(C_{2k+1}, r)$, we need to construct a $C_{2k+1}$-free graph $G$ with the maximum $\delta(G)$ such that $\chi(G)\ge r+1$.   Let $G_{r+1}$ denote the graph obtained by taking $r+1$ vertex-disjoint copies of $K_{\frac{n}{2(r+1)},\frac{n}{2(r+1)}}$ and selecting a vertex in each of them such that these vertices form a $K_{r+1}$. See  Fig.3 below.
\begin{figure}[H]
\centering
\includegraphics[width=8cm]{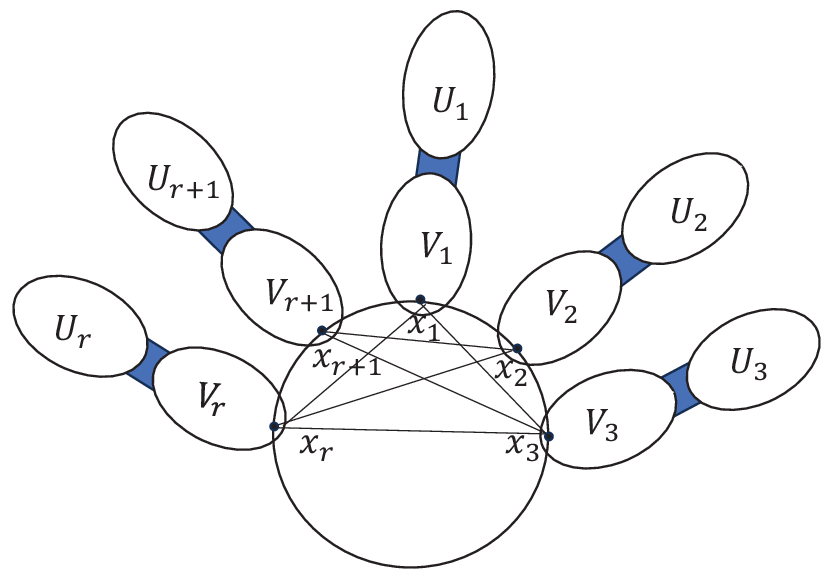}
\caption{$G_{r+1}$}
\end{figure}

Note that $G_{r+1}$ is $C_{2k+1}$-free (if $2k+1>r+1$) and $G_{r+1}$ is not $r$-partite since it contains $K_{r+1}$, so $\delta_{\chi}(C_{2k+1}, r)\ge \delta(G_{r+1})/n=\frac{1}{2r+2}$. We will show that $G_{r+1}$ is the unique extremal construction for $\delta_{\chi}(C_{2k+1}, r)$.  We will prove the following main results implying that $\delta_{\chi}(C_{2k+1}, r)\le \delta(G_{r+1})/n=\frac{1}{2r+2}$ if $r\ge 3$ and $k\ge 3r+4$.
\begin{thm}\label{main}
Let $r, k$ and $n$ be integers with $r\ge 3$,  $k\ge 3r+4$ and $n\ge 108(r+1)^rk$. Let $G$ be an $n$-vertex $C_{2k+1}$-free graph. If
$$\delta(G)\geq \frac{n}{2r+2},$$
then $G$ is r-partite, or $G=G_{r+1}$.
\end{thm}

Theorem \ref{main} implies the following corollaries.

\begin{cor}\label{mainchi}
 $\delta_{\chi}(C_{2k+1}, r)=\frac{1}{2r+2}$ for $r\ge 3$ and $k\ge 3r+4$.
 \end{cor}

\begin{cor}
Let  $r\ge 3$ be an integer.   If $G\neq G_{r+1}$ is  a non $r$-partite $n$-vertex  graph with
$\delta(G)\geq \frac{n}{2r+2}$,
then $C_{2k+1} \subset G$ for all $k\in [3r+4, \frac{1}{108(r+1)^r}n]$.
\end{cor}

Let $r\ge 3$ and ${\mathcal C}$ be a family consisting of some  odd cycles in which $C_{2p+1}$ is the shortest  cycle  in ${\mathcal C}$ and $C_{2k+1}$ is the longest  cycle in ${\mathcal C}$ satisfying $2p\ge r+1$ and $k\ge 3r+4$. Since any ${\mathcal C}$-free graph is $C_{2k+1}$-free,  $\delta_{\chi}({\mathcal C}, r)\le \delta_{\chi}(C_{2k+1}, r)=\frac{1}{2r+2}$. On the other hand, $G_{r+1}$ does not contain any odd cycle longer than $r+1$, and the shortest cycle in ${\mathcal C}$   has length $2p+1> r+1$, so $G_{r+1}$ is ${\mathcal C}$-free as well. Note that $G_{r+1}$ has the minimum degree $\frac{n}{2r+2}$ and $\chi(G_{r+1})=r+1$, thus, $\delta_{\chi}({\mathcal C}, r)\ge \delta(G_{r+1}/n)=\frac{1}{2r+2}$. To summarize, we have the following corollary.

\begin{cor}\label{cor2}
Let $r\ge 3$ and ${\mathcal C}$ be a family of odd cycles in which $C_{2p+1}$ is the shortest  cycle  in ${\mathcal C}$ and $C_{2k+1}$ is the longest  cycle in ${\mathcal C}$ satisfying $2p\ge r+1$ and $k\ge 3r+4$. Then
 $\delta_{\chi}({\mathcal C}, r)=\frac{1}{2r+2}$.
 \end{cor}

 For $r=2$ and any family consisting  of  some odd cycles, we obtain an improvement on the condition of $n$ in Theorem \ref{YuanPeng} when $k\ge 4p+1$.
 \begin{thm}\label{main2}
 Let $k\ge 4p+1$ and $n\ge 108(2p+1)^{2p}k$ be positive integers.  Let ${\mathcal C}$ be a family of odd cycles in which $C_{2p+1}$ is the shortest odd cycle not in ${\mathcal C}$ and $C_{2k+1}$ is the longest odd cycle in ${\mathcal C}$.  
 If $G$ is an $n$-vertex ${\mathcal C}$-free graph with $\delta(G)\ge \frac{n}{2(2p+1)}$, then $G$ is bipartite or $G=BC_{2p+1}(n)$.
 \end{thm}
 Note that $\max\{ \frac{n}{2(2p+1)}, \frac{2}{2k+3}n\}=\frac{n}{2(2p+1)}$ if and only if  $k\ge 4p+1$, thus
 the condition $k\ge 4p+1$ in Theorem \ref{main2}  is tight.

 Throughout the paper, let $P_{xy}$ denote a path with end points $x$ and $y$. We call a path $P$ an {\em even (odd)} path if $|V(P)|$ is even (odd). The innovation and essence of our proof is to introduce the following concept `strong-2k-core' which seems to be perfect for our problem. We think that how to define a proper `core' might be the key to this type of questions.

\begin{defi}
Let $H$ be a subgraph of $G$. We call $H$ a {\em $2k$-core of $G$} if for each pair of vertices $x, y\in V(H)$, there exists an even path $P_{xy}$ in $H$ of order (the number of vertices) at most $2k$.  We call $H$ a {\em strong-$2k$-core of $G$} if for each pair of vertices $x, y\in V(H)$, there exists an even path $P_{xy}$ in $H$ of order at most $2k$ and there exists an odd path $P'_{xy}$ in $H$ of order at most $2k$.
\end{defi}

To obtain $\delta_{\chi}(C_{2k+1}, 2)$ in Theorem \ref{YuanPeng} (\cite{YuanPeng, YuanPeng2}), we assume on the contrary that $\chi (G)\ge 3$ under the condition of Theorem \ref{YuanPeng}. Since $\chi (G)\ge 3$, we can take a shortest cycle $C$ and analyze the structure of the vertices outside $C$ having distances $1$ or $2$ to $C$. This will not work when $c\ge 3$, and it is somehow natural since a shortest cycle will not take the full advantage of the condition $\chi(G)\ge c+1\ge 4$.

To overcome the above barrier, we  come up with the concept  strong-$2k$-core
which is very crucial to give the solution for $\delta_{\chi}(C_{2k+1}, c)$ if $c\ge 3$. A maximum strong $2k$-core seems to be `perfect' for analyzing structures of $C_{2k+1}$-free graphs with chromatic number greater than $3$ and some minimum degree conditions, and  it is the key that we can give a concise proof  while methods used in \cite{YuanPeng, YuanPeng2} and other related papers cannot be applied to obtain Theorems \ref{main} and \ref{main2}. Indeed, we unify the proofs of Theorems \ref{main} and \ref{main2} by taking out a maximum strong-2k-core. We think that this concept grasps the essence of the problem and it makes  our proof  concise and elementary (we do not need to borrow any other tools).
In the same time, we have also applied this  new concept to
obtain strong structural stability result for $C_{2k+1}$-free graphs in
 another paper \cite{cycle-strong} (see more details in Section \ref{remark}).  How to define a proper `core' might be a  key to this type of questions.

We follow standard notations through. For $S\subseteq V(G)$, let $G[S]$ denotes the subgraph of $G$ induced by $S$, and let $G-S$ denote the subgraph induced by $V(G)-S$. For any vertex $v\in V(G)$, let $N(v)$ denote the set of all neighbors of $v$ in $G$, and $N_S(v)=N(v)\cap S$. Let $d(v)=\vert N(v)\vert$ and $d_S(v)=\vert N_S(v)\vert$. 
Let $N_G(S)$ denote the union of  the neighborhoods of all vertices of $S$ in $G$, and we simply write $N(S)$ sometimes. 

\section{Proofs of Theorems \ref{main} and \ref{main2}}
We first  prove some crucial lemmas.

\begin{lem}\label{mainlemma}
Let $r\ge 2$, $k$ and $n$ be integers with $n\ge108(r+1)^rk$ and $k\ge f(r)$, where $f(r)=2r+1$ if $r=2$ and $f(r)=3r+4$ if $r\ge 3$. Let $G$ be a $C_{2k+1}$-free graph on $n$ vertices with $\delta(G)\ge \frac{n}{2r+2}$, then for any even path $P_{xy}$ of order at most $2k$, we have $|(N(x)\cap N(y))\setminus V(P_{xy})|\le 15r.$
\end{lem}
\noindent\emph{\textbf{Proof of Lemma \ref{mainlemma}}.} We show that if there exists an even path $P_{xy}$ of order at most $2k$ satisfying $|(N(x)\cap N(y))\setminus V(P_{xy})|>15r$, then there is a $C_{2k+1}$ in $G$. We prove it by  induction on $2k-|V(P_{xy})|$. If $|V(P_{xy})|=2k$, then it is clear that there is a $C_{2k+1}$ in $G$ if $|(N(x)\cap N(y))\setminus V(P_{xy})|\ge1$. Suppose that it holds for $2k\ge |V(P_{xy})|\ge 2l+2$. Next, we will show that it holds for $|V(P_{xy})|=2l$. Assume that there exist an even path $P_{xy}$ of order $2l$ such that $|(N(x)\cap N(y))\setminus V(P_{xy})|> 15r.$ Let $U\subset (N(x)\cap N(y))\setminus V(P_{xy})$ with $|U|=15r$. If there exists $u, v\in U$ such that $|(N(u)\cap N(v))\setminus V(P_{xy})|> 15r$, then there is an even path $P_{uv}=uP_{xy}v$ of order $2l+2$ such that $|(N(u)\cap N(v))\setminus V(P_{uv})|> 15r$. By induction hypothesis there is a $C_{2k+1}$ in $G$. Therefore we may assume that $|(N(u)\cap N(v))\setminus V(P_{xy})|\le 15r$ for any $u, v\in U$. Combining with $\delta(G)\ge \frac{n}{2r+2}$, we have
\begin{eqnarray*}\label{egl}
\vert V(G)\vert&\ge &\bigg(\frac{n}{2r+2}-|V(P_{xy})|\bigg)|U|-{|U|\choose 2}15r \\
&\ge&\bigg(\frac{n}{2r+2}-2k\bigg)15r-\frac{(15r)^3}{2}\\
&>& n
\end{eqnarray*}
by  $k\ge f(r)$, $n\ge 108(r+1)^rk$ and direct calculation.
This is a contradiction.
$\hfill\square$

\begin{lem}\label{core1}
Let $r\ge 2$, $k$ and $n$ be integers with $n\ge108(r+1)^rk$ and $k\ge f(r)$, where $f(r)=2r+1$ if $r=2$ and $f(r)=3r+4$ if $r\ge 3$.
Let $G$ be a $C_{2k+1}$-free graph on $n$ vertices with $\delta(G)\ge \frac{n}{2r+2}$. If $H$ is a $2k$-core of $G$, then $|V(H)|\le 2r+2$.
\end{lem}
\noindent\emph{\textbf{Proof of Lemma \ref{core1}}.} Suppose that $|V(H)|\ge 2r+3$, then take $U\subset V(H)$ such that $|U|=2r+3$. Since $H$ is a $2k$-core of $G$, by Lemma \ref{mainlemma}, $|(N(x)\cap N(y))\setminus V(P_{xy})|\le 15r$ for $x, y\in U$, where $P_{xy}$ is an even path of order at most $2k$. Therefore, for $r\ge 3$, $|N(x)\cap N(y)|\le 15r+2k-2\le 7k$ since $k\ge 3r+4$. Combining with $\delta(G)\ge \frac{n}{2r+2}$, we have
\begin{eqnarray*}\label{egl}
\vert V(G)\vert&\ge &\frac{n}{2r+2}|U|-{|U|\choose 2}7k \\
&=&\frac{(2r+3)n}{2r+2}-{2r+3\choose 2}7k=n+\frac{n}{2r+2}-{2r+3\choose 2}7k \\
&>&n
\end{eqnarray*}
by  $n\ge 108(r+1)^rk$ and direct calculation.
This is a contradiction.

For $r=2$, we have $|N(x)\cap N(y)|\le 15r+2k-2=2k+28$. Combining with $\delta(G)\ge \frac{n}{2r+2}=\frac{n}{6}$, we have
\begin{eqnarray*}\label{egl}
\vert V(G)\vert&\ge &\frac{n}{6}|U|-{|U|\choose 2}(2k+28) \\
&=&\frac{7n}{6}-{7\choose 2}(2k+28)=n+\frac{n}{6}-(42k+588) \\
&>&n
\end{eqnarray*}
by $n\ge 108(r+1)^rk$, $k\ge 5$ and direct calculation.
This is a contradiction. $\hfill\square$

The following observation is easy to see, but it is  important to our proof.

\begin{fact}\label{max}
Let $H$ be a strong-$2k$-core of $G$ with $|V(H)|=l\le 2k-2$. If there exists an even path $P_{uv}\subset V(G)\setminus V(H)$ such that $xu\in E(G)$ and $xv\in E(G)$ for some $x\in V(H)$ with $|V(P_{uv})|\le 2k-l$, then $V(H)\cup V(P_{uv})$ is a strong-$2k$-core of $G$. If there exists a path $P_{uv}\subset V(G)\setminus V(H)$ (it is possible that $u=v$) such that $xu\in E(G)$ and $yv\in E(G)$ for some $x, y\in V(H)$ with $|V(P_{uv})|\le 2k-l$, then $V(H)\cup V(P_{uv})$ is a strong-$2k$-core of $G$.
\end{fact}

\begin{lem}\label{core2}
Let $r\ge 2$, $k$ and $n$ be integers with $n\ge108(r+1)^rk$ and $k\ge f(r)$, where $f(r)=2r+1$ if $r=2$ and $f(r)=3r+4$ if $r\ge 3$.
Let $G$ be a $C_{2k+1}$-free graph on $n$ vertices with $\delta(G)\ge \frac{n}{2r+2}$. If $H$ is a strong-$2k$-core of $G$, then $|V(H)|\le r+1$.
\end{lem}
\noindent\emph{\textbf{Proof of Lemma \ref{core2}}.} Let $H$ be a maximum strong-$2k$-core of $G$, i.e., $|V(H)|$ has the maximum cardinality among all strong-$2k$-cores of $G$. Suppose on the contrary that $V(H)=\{x_1, x_2,\dots, x_l\}$ with $l\ge r+2$. By Lemma \ref{core1}, we have that $r+2\le l\le 2r+2$.  Let
$$N_i=N(x_i)\setminus V(H) \ \text{for} \ 1\le i\le l.$$ By the maximality of $|V(H)|$ and Fact \ref{max}, we have that
\begin{claim}\label{claim1}
$N_i$ is an independent set for $1\le i\le l$ and $N_i\cap N_j=\emptyset$ for $1\le i<j\le l$.
\end{claim}
Fix $y_i\in N_i$ for $1\le i\le l$ and let $$N_i^*=N(y_i)\setminus V(H).$$
Claim \ref{claim1} says that $N_i$ is an independent set, so we have that
\begin{claim}\label{claim2}
 $N^*_i\cap  N_i =\emptyset$ for $1\le i\le l$.
 \end{claim}
By the maximality of $|V(H)|$ and Fact \ref{max}, we have that
\begin{claim}\label{claim3}
$N^*_i\cap N_j=\emptyset \ \text{and} \ N^*_i\cap N^*_j=\emptyset \ \text{for} \ 1\le i<j\le l.$
\end{claim}
Since $\delta(G)\ge \frac{n}{2r+2}$, we have that $|N_i|\ge \frac{n}{2r+2}-l$. Claim \ref{claim1} says that $N_i\cap N_j=\emptyset$ for $1\le i<j\le l$, this implies that  $|N(y_i)\cap H|=1$, so $|N^*_i|\ge \frac{n}{2r+2}-1$ for $1\le i\le l$. Combining with Claims \ref{claim1}, \ref{claim2}, \ref{claim3}, we have
$$|V(G)|\ge\sum_{i=1}^l(|N_i|+|N^*_i|)\ge (r+2)\bigg(\frac{n}{r+1}-(2r+3)\bigg)>n+\frac{n}{r+1}-(r+2)(2r+3)>n$$
since  $n\ge 108(r+1)^rk$ and $k\ge f(r)$.
A contradiction. $\hfill\square$

Note that an odd cycle with length less than $2k+1$ in a graph $G$  is a strong-$2k$-core of $G$. Next we show that there exists an odd cycle with length less than $2k+1$ in a graph $G$ satisfying the conditions in Theorem \ref{main} or Theorem \ref{main2}.
\begin{lem}\label{shortcycle}
Let $r\ge 2$, $k$ and  $n$ be positive integers.
If  $G$ is a  $C_{2k+1}$-free non-bipartite graph on $n$ vertices with $\delta(G)\ge \frac{n}{2r+2}$, then $G$ contains an odd cycle with length no more than $2(2r+1)+1$.
\end{lem}
\noindent\emph{\textbf{Proof of Lemma \ref{shortcycle}}.}
Let $C_{2m+1}=v_1v_2\cdots v_{2m+1}v_1$ be a shortest odd cycle of $G$. Let
$$G'=G-V(C_{2m+1}).$$

\begin{claim}\label{claim21}
For any vertex $v\in V(G')$, we have $d_{C_{2m+1}}(v)\leq 2 \mbox{~if~} m\geq 2.$
\end{claim}
\noindent\emph{\textbf{Proof of Claim \ref{claim21}}.} Let $m\geq 2$.  Suppose on the contrary that there exists a vertex $x\in V(G')$, such that $d_{C_{2m+1}}(x)\geq 3$, let $\{v_i,v_j,v_q\}\subseteq N_{C_{2m+1}}(x)$, where $1\leq i<j<q\leq 2m+1$. We claim that any two vertices of $\{v_i,v_j,v_q\}$ are not adjacent. Otherwise, without loss of generality,  assume that $v_iv_j\in E(G)$, then $vv_iv_jv$ is a copy of $C_3$, a contradiction to $C_{2m+1}$ being a shortest cycle. Moreover, $C_{2m+1}$ is divided into three paths by $\{v_i,v_j,v_q\}$, since $C_{2m+1}$ is an odd cycle of $G$, there is at least one even path (so the length is odd). Without loss of generality, assume that $v_qv_{q+1}\cdots v_{2m+1}v_1\cdots v_i$ is an even path of $C_{2m+1}$. We have shown that any two vertices of $\{v_i,v_j,v_q\}$ are not adjacent, so $v_iv_{i+1}\cdots v_{j}v_{j+1}\cdots v_q$ is an odd path with at least $5$ vertices, then we use the odd path $v_ivv_q$ to replace the odd path $v_iv_{i+1}\cdots v_{j}v_{j+1}\cdots v_q$ of $C_{2m+1}$ to get a shorter odd cycle $v_ivv_qv_{q+1}\cdots v_{2m+1}v_1\cdots v_i$, a contradiction. This completes the proof of Claim \ref{claim21}. $\hfill\square$

Since $C_{2m+1}$ is a shortest odd cycle of $G$, it does not contain any chord. By Claim \ref{claim21} and $\delta(G)\geq \frac{n}{2r+2}$, for $m\geq 2$,
$$(2m+1)\cdot (\frac{n}{2r+2}-2)\leq e(C_{2m+1},G')\leq 2\cdot (n-2m-1).$$
This implies that $m< 2r+2$. Since $m$ is an integer, $m\le 2r+1$. This completes the proof of Lemma \ref{shortcycle}. $\hfill\square$

If $G$ is a  $C_{2k+1}$-free non-bipartite graph on $n$ vertices with $\delta(G)\ge \frac{n}{2r+2}$ and $k\ge 2r+1$, then by Lemma \ref{shortcycle}, $G$ contains an odd cycle with length no more than $2k-1$. Therefore, $G$ contains a strong-$2k$-core with at least $3$ vertices. Combining with Lemma \ref{core2}, we have the following corollary.

\begin{cor}\label{existcore}
 Let $r\ge 2$, $k$ and  $n$ be positive integers with $k\ge 2r+1$.
If  $G$ is a  $C_{2k+1}$-free non-bipartite graph on $n$ vertices with $\delta(G)\ge \frac{n}{2r+2}$, then $G$ contains a strong-$2k$-core $H$ with $3\le \vert V(H)\vert \le r+1$.
\end{cor}

Now we are ready to give structural information on graphs satisfying the conditions in Theorem \ref{main} or Theorem \ref{main2}.

\begin{lem}\label{structure}
Let $r\ge 2$, $k$ and $n$ be integers with $n\ge108(r+1)^rk$ and $k\ge f(r)$, where $f(r)=2r+1$ if $r=2$ and $f(r)=3r+4$ if $r\ge 3$.
Let $G$ be a non-bipartite $C_{2k+1}$-free graph on $n$ vertices with $\delta(G)\ge \frac{n}{2r+2}$. Let $H$ be a maximum strong-$2k$-core of $G$ with  $V(H)=\{x_1, x_2,\dots, x_l\}$ and $3\le l\le r+1$.  For each $i$,  $1\le i\le l$, let
$$N_i=N(x_i)\setminus V(H),  \ N^1_i=N(N_i)\ \text{and} \ N^2_i=N(N^1_i)\setminus V(H).$$
Then the following holds.

(i) If  $l=r+1$, then for each $1\le i\le l$, $N^1_i$ and $N^2_i$ are independent sets, $|N^1_i|=\frac{n}{2r+2}$, $|N^2_i|=\frac{n}{2r+2}$ and $G[N^1_i, N^2_i]$ forms a complete bipartite graph.

(ii) If  $l\le r$, then $x_i$ is a cut vertex for each $1\le i\le l$.
\end{lem}

\noindent\emph{\textbf{Proof of Lemma \ref{structure}}.}
By the maximality of $|V(H)|$ and Fact \ref{max}, we have that
\begin{equation}\label{eq1}
N_i \ \text{is an independent set and} \ N_i\cap N_j=\emptyset \ \text{for} \ 1\le i<j\le l.
\end{equation}
Since $N_i$ is an independent set, we have that
\begin{equation}\label{eq2}
 N^1_i\cap N_i=\emptyset \ \  \text{for} \ 1\le i\le l.
 \end{equation}
By the maximality of $|V(H)|$ and Fact \ref{max}, we have that
\begin{equation}\label{eq3}
N^1_i\cap N_j=\emptyset, \ N^1_i\cap N^1_j=\emptyset,  \ N_i\cap N^2_j=\emptyset, \  \text{and} \ N^1_i\cap N^2_j=\emptyset \ \text{for} \ 1\le i\neq j\le l.
\end{equation}

Since $\delta(G)\ge \frac{n}{2r+2}$, we have
\begin{equation}\label{eq6}
|N^1_i|\ge \frac{n}{2r+2} \ \   \text{for}  \ 1\le i\le l.
\end{equation}
By the maximality of $|V(H)|$ and Fact \ref{max},  for any vertex $y \in N^1_i$ and $y\neq x_i$, we have that $N(y)\cap V(H)=\emptyset$. Combining with   $\delta(G)\ge \frac{n}{2r+2}$, we have
\begin{equation}\label{eq7}
|N^2_i|\ge \frac{n}{2r+2} \ \   \text{for}  \ 1\le i\le l.
\end{equation}

\textbf{Case (i).} $l=r+1$.

\begin{claim}\label{claimindependent}
$N^1_i$ is an independent set for $1\le i\le l$.
\end{claim}
\noindent\emph{\textbf{Proof of Claim \ref{claimindependent}}.} Without loss of generality, assume that $N^1_1$ is not an independent set, and $u_1u_2\in E(G)$ for some $u_1, u_2\in N^1_1$. Let $u_3, u_4\in N_1$ such that $u_1u_3, u_2u_4\in E(G)$. By Fact \ref{max}, we have $u_3=u_4$. So $\{u_1, u_2, u_3\}$ forms a $K_3$. By Lemma \ref{mainlemma}, $|N(u_i)\cap N(u_{i+1})|\le 15r$ for $i\in [3]$. Hence $|N_1^1\cup N_1^2|\ge |N(u_1)\cup N(u_2)\cup N(u_3)| \ge\frac{3n}{2r+2}-45r$. Combining with (\ref{eq1}), (\ref{eq2}), and (\ref{eq3}),  we have
\begin{eqnarray*}\label{egl}
\vert V(G)\vert&\ge &\sum_{i=2}^{l}(|N_i|+|N^1_i|)+|N_1^1\cup N_1^2|\\
&\ge&r\bigg(\frac{n}{r+1}-r\bigg)+\frac{3n}{2r+2}-45r \\
&>&n-\frac{n}{r+1}-r^2+\frac{3n}{2r+2}-45r \\
&>&n
\end{eqnarray*}
since $k\ge f(r)$ and $n\ge 108(r+1)^rk$.
A contradiction.
$\hfill\square$

By Claim \ref{claimindependent}, we have
\begin{equation}\label{eq4}
 N^1_i\cap N^2_i=\emptyset \ \text{for} \ 1\le i\le l.
\end{equation}
 Combining (\ref{eq4}), (\ref{eq3}),   (\ref{eq6}), and  (\ref{eq7}),
we have
 $$n=|V(G)|\ge \sum_{i=1}^{l}(|N^1_i|+|N^2_i|)\ge(r+1)\frac{2n}{2r+2}=n.$$
 Therefore,
$|N^1_i|=\frac{n}{2r+2}$, $|N^2_i|=\frac{n}{2r+2}$ and $G[N^1_i, N^2_i]$ forms a complete bipartite graph for $1\le i\le l$. Since $G$ is $C_{2k+1}$-free, $N^2_i$ is an independent set.

\textbf{Case (ii).} $l\le r$.

Since $l\ge 3$, in this case, $r\ge 3$ and $k\ge 3r+4$.

We show that $x_i$ is a cut vertex for each $1\le i\le l$. Otherwise,
without loss of generality, suppose that there exists a path $P_{uv}\subseteq V(G)\setminus V(H)$ such that $x_1u\in E(G)$ and $x_2v\in E(G)$ for $x_1, x_2\in V(H)$, and assume that $|V(P_{uv})|$ has the minimum cardinality. Let $A=V(H)\cup V(P_{uv})$. By the maximality of $|V(H)|$ and Fact \ref{max}, we have $|A|\ge 2k+1$.
By the minimality of $|V(P_{uv})|$ and $|V(H)|\le r$, then $d_{A}(x)\le r$ for $x\in A$. By the minimality of $|V(P_{uv})|$ and the maximality of $|V(H)|$ and Fact \ref{max}, we have that $d_{A}(x)\le 3$ for $x\in V(G)\setminus A$. Hence
$$|A|\cdot (\frac{n}{2r+2}-r)\leq e(A,V(G)\setminus A)\leq 3(n-|A|),$$
then $|A|\le\frac{(6r+6)n}{n-2r(r+1)+3}$. Recall that $|A|\ge 2k+1$, $k\ge 3r+4$ and
$n\ge 108r^{r-1}k$, then we get  $n<0$, a contradiction. Therefore $x_i$ is a cut vertex for $1\le i\le l$.
$\hfill\square$

\bigskip
Now we are ready to prove Theorems \ref{main} and \ref{main2}.

\bigskip

\noindent\emph{\textbf{Proof of Theorem \ref{main2}}.} Let $k\ge 4p+1$ and $n\ge 108(2p+1)^{2p}k$ be positive integers.
Let ${\mathcal C}$ be a family of odd cycles in which $C_{2p+1}$ is the shortest odd cycle not in ${\mathcal C}$ and $C_{2k+1}$ is the longest odd cycle in ${\mathcal C}$.  Recall that $BC_{2p+1}(n)$ is the graph obtained by taking $2p+1$ vertex-disjoint copies of $K_{\frac{n}{2(2p+1)},\frac{n}{2(2l+1)}}$ and selecting a vertex in each of them such that these vertices form a cycle of length $2p+1$.  Let $G$ be an $n$-vertex ${\mathcal C}$-free non-bipartite graph with $\delta(G)\ge \frac{n}{2(2p+1)}$. We are going to show that $G=BC_{2p+1}(n)$.

Take $r=2p$ in Lemma \ref{structure}.  Notations used in the proof  follow from Lemma \ref{structure}.  For example, $H$ is a maximum strong-$2k$-core of $G$ with  $V(H)=\{x_1, x_2,\dots, x_l\}$ and $3\le l\le 2p+1$.  We claim that $l=2p+1$. Otherwise,  $l<2p+1$, since $H$ is a maximum strong-$2k$-core of $G$, there must be an odd cycle in $H$ and the length of this cycle is no more than $l<2p+1$. Recall that $C_{2p+1}$ is the shortest odd cycle not in ${\mathcal C}$, in other words, ${\mathcal C}$ contains all odd cycles with length less than $2p+1$. Thus, $H$ is not ${\mathcal C}$-free, a contradiction. So we have $l=2p+1$.

Applying Lemma \ref{structure} (i), we obtain that for each $1\le i\le l$, $N^1_i$ and $N^2_i$ are independent sets, $|N^1_i|=\frac{n}{2r+2}$, $|N^2_i|=\frac{n}{2r+2}$ and $G[N^1_i, N^2_i]$ forms a complete bipartite graph.  Since $H$ is a strong-$2k$-core of $G$ with $2p+1$ vertices and $H$ does not contain any odd cycle shorter than $2p+1$, $H$ must be $C_{2p+1}$. Therefore $G=BC_{2p+1}(n)$. This completes the proof of Theorem \ref{main2}. $\hfill\square$

\noindent\emph{\textbf{Proof of Theorem \ref{main}}.}  Let $r\ge 2$, $k$ and $n$ be integers with $n\ge108(r+1)^rk$ and $k\ge f(r)$, where $f(r)=2r+1$ if $r=2$ and $f(r)=3r+4$ if $r\ge 3$. Let $G$ be a $C_{2k+1}$-free graph on $n$ vertices with $\delta(G)\ge \frac{n}{2r+2}$, we show that $G$ is $r$-partite, or $G=G_{r+1}$.

We apply induction on $r$. When $r=2$, applying Theorem \ref{main2} with ${\mathcal C}=C_{2k+1}$ and $2p=2$, we are fine. We will show that the conclusion holds for $r\ge 3$. We will apply Lemma \ref{structure} and  notations used in the proof  follow from Lemma \ref{structure}. For example, $H$ is a maximum strong-$2k$-core of $G$ with  $V(H)=\{x_1, x_2,\dots, x_l\}$ and $3\le l\le r+1$.

\textbf{Case (i).} $l=r+1$.

By Lemma \ref{structure},   for each $1\le i\le l$, $N^1_i$ and $N^2_i$ are independent sets, $|N^1_i|=\frac{n}{2r+2}$, $|N^2_i|=\frac{n}{2r+2}$ and $G[N^1_i, N^2_i]$ forms a complete bipartite graph.   If $H=K_{r+1}$, then $G=G_{r+1}$. If $H\neq K_{r+1}$, then $\chi(G)\le r$.

\textbf{Case (ii).} $l\le r$.

By Lemma \ref{structure}, $x_i$ is a cut vertex for $1\le i\le l$. Let
$$H_i=\{u\in V(G)\setminus V(H)| \ \text{there exists a path} \ P_{ux_i}\setminus\{x_i\}\subseteq V(G)\setminus V(H)\} \ \text{for} \ 1\le i\le l.$$
Clearly,  $V(G)$ is partitioned into $V(G)=(\cup_{i=1}^l H_i)\cup V(H)$. Therefore it is sufficient to show that $G[H_i]$ is ($r-1$)-partite. Note that $N_i\cup N_i^1\subset H_i\cup \{x_i\}$, by (\ref{eq2}), $|H_i|\ge \frac{2n}{2r+2}-l\ge 108r^{r-1}k$. Note that $l\ge 3$, then $|H_i|\le n-\frac{4n}{2r+2}+l=\frac{(r-1)n}{r+1}+l$. Note that $N(x)\subset H_i\cup\{x_i\}$ for any $x\in H_i$, hence $\delta (G[H_i])\ge\delta(G)-1=\frac{n}{2r+2}-1> \frac{|H_i|}{2r}$. By induction hypothesis, $G[H_i]$ is ($r-1$)-partite. This completes the proof.
$\hfill\square$

\section{Remarks}\label{remark}

To determine $\delta_{\chi}(\mathcal{F}, c)$, we need to find the maximum $\delta(G)$ an $\mathcal{F}$-free graph $G$ with large enough number of vertices and $\chi(G)>c$ can have. Will our construction $G_{r+1}$ in Section 2 (see Figure 3)  be the extremal construction for more  values of $k$ and  $r\ge 3$? It would be interesting to improve the condition on $k$ and $n$ in Theorem \ref{main}.

 The innovation of this paper is to give the concept  strong-$2k$-core,  it is the key that we can give a concise proof  while other previous methods  cannot be applied to obtain our results.  Taking a maximum  strong-$2k$-core seems to be `perfect' for analyzing structures of $C_{2k+1}$-free graphs for our question.  We have obtained another application of this concept to give a strong structural stability result for $C_{2k+1}$-free graphs in \cite{cycle-strong}.  F\"uredi and Gunderson\cite{FuGun} showed that  $ex(n, C_{2k+1})$ is achieved only on $K_{\lfloor\frac{n}{2}\rfloor, \lceil\frac{n}{2}\rceil}$ if $n\ge 4k-2$. It is natural to study how far a $ C_{2k+1}$-free graph is from being bipartite. If a graph $G$ and a graph $H$ have at most one vertex in common and there is no edge connecting $V(G)-V(G)\cap V(H)$ and $V(H)-V(G)\cap V(H)$, then we call graph $H$ a {\em suspension} to graph $G$ with $\vert V(G)\cap V(H)\vert$ suspension point. Let $T^*(r, n)$  be obtained by adding a suspension $K_{r}$ with $1$ suspension point to $K_{\lfloor\frac{n-r+1}{2}\rfloor, \lceil\frac{n-r+1}{2}\rceil}$. By taking a maximum  strong-$2k$-core of a graph, we showed the following result in \cite{cycle-strong}.
 \begin{thm}\label{main-stability}\cite{cycle-strong}
Let $r\ge 1$,  $2k\ge r+4$ and $n\ge 20(r+2)^2k$ be integers. Let $G$ be an $n$-vertex $C_{2k+1}$-free graph. If $e(G)\ge \big\lfloor{\frac{(n-r+1)^2}{4}}\big\rfloor+{r\choose 2}= e(T^*(r, n))$, then $G$ is obtained by adding suspensions to a bipartite graph $B=G[V_1, V_2]$ one by one (in other words, $G=B\bigcup\limits_{i=1}^p G_i$ for some $p$,  $G_1$ is a suspension to $B$,  $G_j$ is a suspension to $B\bigcup\limits_{i=1}^{j-1} G_i$ for $2\le j\le p$)  such that the total number of vertices not in $B$ is no  more than $r-1$. Furthermore, the total number of vertices not in $B$ equals $r-1$ if and only if $G=T^*(r, n)$.
\end{thm}
 Let $$d_2(G)=\min\{|T|: T\subseteq V(G), G-T \ \text{is bipartite}.\},$$
  $$\gamma_2(G)=\min\{|E|: E\subseteq E(G), G-E \ \text{is bipartite}.\}.$$
 Theorem \ref{main-stability}  implies that $d_2(G)\le r-1$ (removing all vertices not in $B$ in Theorem \ref{main-stability} yields a bipartite graph),  and $\gamma_2(G)\le {\lceil\frac{r}{2}\rceil \choose 2}+{\lfloor\frac{r}{2}\rfloor \choose 2}$ (the graph induced by the set of vertices not in $B$ in Theorem \ref{main-stability} has a bipartite graph with at least half of the edges), which is a recent result of Ren-Wang-Wang-Yang\cite{RWWY}. In \cite{RWWY}, they proved $d_2(G)\le r-1$ and $\gamma_2(G)\le {\lceil\frac{r}{2}\rceil \choose 2}+{\lfloor\frac{r}{2}\rfloor \choose 2}$ in two separate theorems by different proofs. By taking a maximum  strong-$2k$-core, we \cite{cycle-strong} give a new and simpler method to obtain Theorem \ref{main-stability}  implying both results in \cite{RWWY}. It would be  interesting to understand what should be taken out for many unsolved questions. We think that defining a proper `core' might be a key to this type of questions.

A related concept chromatic threshold concerns how large the minimum degree needs  to guarantee that the chromatic number of $\mathcal F$-free graphs is bounded by some constant.  Precisely,  the chromatic threshold of $\mathcal F$ is defined by
$$\delta_{\chi}(\mathcal{F})=\inf\{\alpha| \exists \ K \ \text{such that} \ \text{any} \ G\in\mathcal G(n, \mathcal F, \alpha) \ \text{must satisfy}  \ \chi(G)\le K.\}.$$
The chromatic threshold is by now much better understood than the chromatic profile. Building on the work of \L uczak and Thomass\'e \cite{LuTh}, and generalising various previous results, Allen, B\"ottcher, Griffiths, Kohayakawa and Morris \cite{ABGKM} determined the chromatic threshold of every finite family $\mathcal F$. For more details about the history of the study of the chromatic threshold see \cite{ABGKM} and the references therein.

One can consider the more restrictive notion of the homomorphism threshold $\delta_{hom}(\mathcal{F})$ of a family $\mathcal F$, which is  the smallest minimum degree that guarantees that $\mathcal F$-free graphs are homomorphic to a small $\mathcal F$-free graph. That is,
$$\delta_{hom}(\mathcal{F})=\inf\{\alpha| \ \exists \ \text{a graph} \ H \ \text{such that} \ \text{any} \ G\in \mathcal G(n, \mathcal F, \alpha) \ \text{must be homomorphic to} \ H.\}.$$
Note that $\delta_{hom}(\mathcal{F})\ge \delta_{\chi}(\mathcal{F})$. Determining chromatic profile or homomorphism thresholds is distinctively
harder than determining chromatic thresholds. 

\end{document}